# An inequality for anti-self-polar polytopes

Mikhail G. Katz

**Abstract:** We prove an inequality for the $f$-vectors of anti-self-polar polytopes conjectured by Katz in 1989. The proof uses Kalai's combinatorial inequality based on a result of Whiteley. The inequality can also be obtained from the results of Stanley and Karu which however involve difficult algebraic geometry.

We prove an inequality for the $f$-vectors of anti-self-polar polytopes conjectured in (Katz 1989). The proof uses Kalai's combinatorial inequality based on a result of Whiteley (1984). The inequality can also be obtained from the results of Stanley (1987) and Karu (2004) which however involve difficult algebraic geometry.

An anti-self-polar polytope is a polytope $P$ inscribed in the unit sphere, satisfying the relation $P^* = -cP$ for a suitable $c > 0$, where $P^*$ is its polar with respect to the unit sphere. Lovász (1983) exploited anti-self-polar polytopes to prove a combinatorial conjecture. The article Katz (1989) investigated 4-dimensional anti-self-polar polytopes as a possible source of counterexamples to Borsuk's conjecture in $\mathbb{R}^4$ (no counterexamples in this dimension have been found as yet, however).

For an anti-self-polar polytope $P \subseteq \mathbb{R}^4$, we prove the following lower bound on the number of edges in the diameter graph $G(P)$ (i.e., the graph whose edges correspond to pairs of vertices of $P$ at maximal distance). Such a bound was originally conjectured in (Katz 1989). The number of $i$-dimensional faces of $P$ is denoted $f_i$. The number of pairs $x \subseteq y$ where $x$ is an $i$-dimensional face and $y$ is a $j$-dimensional face, is denoted $f_{ij}$. A *facet* of a $d$-dimensional polytope is a $(d-1)$-dimensonal face.

**Theorem 1**. *Let $P \subseteq \mathbb{R}^4$ be a 4-dimensional anti-self-polar polytope. Then the number of edges $e(G(P))$ in the diameter graph $G(P)$ is at least $3f_0(P) - 5$.*

*Proof.* Let $P$ be an anti-self-polar polytope of dimension $d + 1$. Let $e(G(P))$ be the number of edges in the graph $G(P)$. Then $f_{0d}(P) = 2e(G(P))$. Indeed, if $(v,w)$ is a pair of vertices of $P$ at maximal distance, then $w$ is a vertex of the facet opposite $v$. Therefore, each such pair is counted twice in $f_{0d}(P)$.

Thus we need to show that

$$f_{03}(P) \geq 6f_0(P) - 10.$$

We will prove such a bound using a result of Kalai's.

**Definition 2**. *Let P be a 4-dimensonal polytope. For a facet $\phi \subseteq P$, $a_\phi^j$ denotes the number of j-gons occurring as faces of the facet $\phi$.*

Let $a^j$ denote the total number of $j$-gons occurring as faces of the polytope $P$. Kalai (1987, Section 12; see also Kalai 1994, Section 4.3) proved that every 4-dimensional polytope satisfies

$$a^4 + 2a^5 + \cdots \geq 4f_0(P) - f_1(P) - 10$$

The gap between the LHS and RHS of this inequality is referred to $g_2(P)$.

Let $\phi$ run through all the facets of $P$. Applying Euler's formula to each facet and using Kalai's bound, we obtain

$$\begin{aligned} f_{03}(P) &= \sum_\phi f_0(\phi) \\ &= \sum_\phi 2 + f_1(\phi) - f_2(\phi) \\ &= \sum_\phi 2 + \frac{1}{2}\left(3a_\phi^3 + 4a_\phi^4 + 5a_\phi^5 + \cdots\right) - f_2(\phi) \\ &= \sum_\phi 2 + \frac{1}{2}f_2(\phi) + \frac{1}{2}\left(a_\phi^4 + 2a_\phi^5 + \cdots\right) \\ &= 2f_3(P) + f_2(P) + \frac{1}{2}\sum_\phi \left(a_\phi^4 + 2a_\phi^5 + \cdots\right) \\ &= 2f_3(P) + f_2(P) + (a^4 + 2a^5 + \cdots) \\ &\geq 2f_3(P) + f_2(P) + 4f_0(P) - f_1(P) - 10 \\ &= \left(2f_3(P) + 4f_0(P)\right) + \left(f_2(P) - f_1(P)\right) - 10 \\ &= 3f_3(P) + 3f_0(P) - 10 \end{aligned}$$

since the Euler characteristic of the 3-sphere vanishes. In the anti-self-polar case, we have $f_0 = f_3$ and therefore one obtains $f_{03} \geq 6f_0(P) - 10$, as required. $\square$

The proof gives that $g_2(P) = 3f_3(P) + 3f_0(P) - 10$ and hence for every 4-polytope, $g_2(P) = g_2(P^*)$.

Some related calculations appear in Bayer (1987). Ziegler (2007, p. 50) notes that Stanley (1987, Corollary 3.2) proved the inequality $f_{03} \geq 3f_0 + 3f_3 - 10$ for 4-polytopes. The paper mentions that the result follows also from the hard Lefschetz theorem for intersection homology. Stanley's proof works for rational polytopes but the hypothesis was removed by Karu (2004); see also Tay et al. (1995).

Qingsong Wang calculated hundreds of examples of anti-self-polar polytopes with spherical diameter $d$ in the range $\arccos(-\frac{1}{4}) < d < \arccos(-\frac{1}{3})$. The calculations were done in PYTHON using the downward gradient flow of the diameter functional.[1] All the examples generated satisfy the boundary case of equality in the inequality of Theorem 1.

**Funding Acknowledgment.** The author was partially supported by Israel Science Foundation grant 743/22.

There are no competing interests.

**Department of Mathematics, Bar Ilan University, Ramat Gan 5290002 Israel**

**katzmik@math.biu.ac.il**

---

[1] See https://ndag.github.io/anti-self-dual-polyhedra/table.html.